\documentclass[12pt]{amsart}
\usepackage{graphicx}
\begin{document}
\newtheorem{theorem}{Theorem}
\newtheorem{prop}[theorem]{Proposition}
\newtheorem{lemma}[theorem]{Lemma}
\newtheorem{claim}[theorem]{Claim}
\newtheorem{cor}[theorem]{Corollary}
\newtheorem{defin}[theorem]{Definition}
\newtheorem{defins}[theorem]{Definitions}
\newtheorem{example}[theorem]{Example}
\newtheorem{xca}[theorem]{Exercise}
\newcommand{\bbb}{\mbox{$\beta$}}
\newcommand{\aaa}{\mbox{$\alpha$}}
\newcommand{\Aaa}{\mbox{$\mathcal A$}}
\newcommand{\ccc}{\mbox{$\mathcal C$}}
\newcommand{\ddd}{\mbox{$\delta$}}
\newcommand{\Ddd}{\mbox{$\Delta$}}
\newcommand{\eee}{\mbox{$\epsilon$}}
\newcommand{\Fff}{\mbox{$\mathcal F$}}
\newcommand{\Ggg}{\mbox{$\Gamma$}}
\newcommand{\ggg}{\mbox{$\gamma$}}
\newcommand{\kkk}{\mbox{$\kappa$}}
\newcommand{\lll}{\mbox{$\lambda$}}
\newcommand{\Lll}{\mbox{$\Lambda$}}
\newcommand{\rrr}{\mbox{$\rho$}}
\newcommand{\sss}{\mbox{$\sigma$}}
\newcommand{\Sss}{\mbox{$\mathcal S$}}
\newcommand{\Ss}{\mbox{$\Sigma$}}
\newcommand{\Th}{\mbox{$\Theta$}}
\newcommand{\ttt}{\mbox{$\tau$}}
\newcommand{\bdd}{\mbox{$\partial$}}
\newcommand{\zzz}{\mbox{$\zeta$}}
\newcommand{\inter}{\mbox{${\rm int}$}}

\title[] {Surfaces, submanifolds, and aligned Fox reimbedding in
non-Haken $3$-manifolds}

\author{Martin Scharlemann}
\address{\hskip-\parindent
         Mathematics Department\\
         University of California\\
         Santa Barbara, CA 93106 \\
         USA}
\email{mgscharl@math.ucsb.edu}

\author{Abigail Thompson}
\address{\hskip-\parindent
         Mathematics Department\\
         University of California\\
         Davis, CA 95616\\
         USA}
\email{thompson@math.ucdavis.edu}

\date{\today}
\thanks{Research supported in part by NSF grants.}

\begin{abstract}
Understanding non-Haken $3$-manifolds is central to many current 
endeavors in $3$-manifold topology.   We describe some results for 
closed orientable surfaces in non-Haken manifolds, and extend Fox's 
theorem for submanifolds of the 3-sphere to submanifolds of general 
non-Haken manifolds.  In the case where the submanifold has connected 
boundary, we show also that the $\bdd$-connected sum decomposition of 
the submanifold can be aligned with such a structure on the 
submanifold's complement.

\end{abstract}

\maketitle

\section{Introduction}

A closed orientable irreducible 3-manifold $N$ is called {\em Haken}
if it contains a closed orientable incompressible surface; otherwise
$N$ is {\em non-Haken}.  In Section 2 we describe some results for
surfaces in non-Haken manifolds.  Generalizing a theorem of Fox
(\cite{F}), we show in Section 3 that a $3$-dimensional submanifold
of a non-Haken manifold $N$ is homeomorphic either to a handlebody
complement in $N$ or the complement of a handlebody in $S^3$. 
Sections 2 and 3 are independent, but both represent progress towards
understanding submanifolds of non-Haken manifolds.  In Section 4 we
combine the techniques from Section 2 with the results from Section 3
to show that if the submanifold $M \subset N$ is $\bdd$-reducible and
has connected boundary, then the embedding can be chosen to align a
full collection of separating $\bdd$-reducing disks in $M$ with
similar disks in the complement of $M$.

\section{Handlebodies in non-Haken manifolds}

Let $N$ be a closed orientable 3-manifold, $F$ a closed orientable
surface of non-trivial genus imbedded in $N$.  Recall that $F$ is {\em
compressible} if there exists an essential simple closed curve on $F$
which bounds an imbedded disk $D$ in $N$ with interior disjoint from
$F$.  $D$ is a {\em compressing disk} for $F$.  

\begin{defin} Suppose $F$ is a separating closed surface in an
orientable irreducible closed $3$-manifold $N$.  $F$ is {\em
reducible} if there exists an essential simple closed curve on $F$
which bounds compressing disks on both sides of $F$.  The union of the
two compressing disks is a {\em reducing sphere} for $F$.

Suppose ${\bf S}$ is a collection of disjoint reducing spheres for
$F$.  A reducing sphere $S \in {\bf S}$
is {\em redundant} if a component of $F - {\bf S}$ that is
adjacent to $S \cap F$ is planar.  ${\bf S}$ is {\em complete} if, for
any disjoint reducing sphere $S'$, $S'$ is redundant in ${\bf S} \cup
S'$.

Let $\sss({\bf S})$ denote the number of components of $F - {\bf
S}$ that are not planar surfaces.  
\end{defin}

Since $N$ is irreducible, any sphere in $N$ is necessarily separating. 
Suppose a reducing sphere $S'$ is added to a collection ${\bf S}$ of
disjoint reducing spheres.  If $S'$ is redundant, the number of
non-planar complementary components in $F$ is unchanged, since $S'$
necessarily separates the component of $F - {\bf S}$ that it
intersects and the union of two planar surfaces along a single
boundary component is still planar.  If $S'$ is not redundant then the
number of non-planar complementary components in $F$ increases by one. 
Thus we have:

\begin{lemma} \label{lemma:complete}
    Suppose ${\bf S} \subset {\bf S'}$ are two collections
of disjoint reducing spheres for $F$ in $N$.  Then $\sss({\bf S}) \leq
\sss({\bf S'})$.  Equality holds if and only if each sphere $S'$ in
${\bf S'} - {\bf S}$ is redundant in $S' \cup {\bf S}$.  In
particular, ${\bf S}$ is complete if and only if for every collection
$ {\bf S'}$ such that ${\bf S} \subset {\bf S'}$, $\sss({\bf S}) =
\sss({\bf S'})$.  \end{lemma} \qed

Let $H$ be a handlebody imbedded in $N$.  $H$ has an {\em unknotted
core} if there exists a pair of transverse simple closed curves $c, d
\subset \partial H$ such that $c \cap d$ is a single point, $d$ bounds
an embedded disk in $H$ and $c$ (the {\em core}) bounds an imbedded
disk in $N$.

\begin {lemma} \label{lemma:pairodisks} Let $F$ be a connected,
closed, separating, orientable surface in a closed orientable
irreducible 3-manifold $N$.  Suppose that $F$ has compressing disks to
both sides.  Then at least one of the following must hold:

\begin{enumerate}
 \item $F$ is a Heegaard surface for $N$.
 \item $N$ is Haken.
 \item There exist disjoint compressing disks for $F$ on opposite sides of $F$.
 \end{enumerate}
\end {lemma}

\begin{proof}
The proof is an application of the generalized Heegaard decomposition 
described in \cite{ST}.    As $F$ is compressible to both sides, we 
can construct a handle decomposition of $N$ starting at $F$ so that 
$F$ appears as a ``thick'' surface in the decomposition.  If $F$ is 
not a Heegaard surface, then this decompostion contains a ``thin'' 
surface $G$ adjacent to $F$. If $G$ is incompressible in $N$, then 
$N$ is Haken.  If $G$ is compressible we apply \cite{CG} to obtain 
the required disjoint compressing disks for $F$.
\end{proof}

\begin{theorem} \label{handlebody}

Let $H$ be a handlebody of genus $g$ imbedded in a closed orientable 
irreducible non-Haken  3-manifold $N$.  Let $G$ be the complement of 
$H$ in $N$.  Let $F=\bdd H=\bdd G$. Suppose $F$ is 
compressible in $G$.  Then at least one of the following must hold:

\begin{enumerate}
 \item The Heegaard genus of $N$ is less than or equal to $g$.
 \item $F$ is reducible.
 \item $H$ has an unknotted core.
 \end{enumerate}

\end{theorem}

\begin{proof}

The proof is by induction on the genus of $H$.  If $g=1$, then the
result of compressing $F$ into $G$ is a $2$-sphere, necessarily
bounding a ball in $N$.  If a ball it bounds lies in $G$ then the
Heegaard genus of $N$ is $\leq 1$.  If a ball it bounds contains
$H$ then $H$ is an unknotted solid torus in $N$ and so it has an 
unknotted core.

Suppose then that $genus(H) = g > 1$ and assume inductively that the
theorem is true for handlebodies of genus $g-1$.  Suppose that $G$,
the complement of $H$, has compressible boundary.  If $G$ is a
handlebody then $G \cup _{F} H$ is a Heegaard splitting of genus $g$
and we are done.  So suppose $G$ is not a handlebody.  Then by Lemma
\ref{lemma:pairodisks} there are disjoint compressing disks on
opposite sides of $F$, say $D$ in $H$ and $E$ in $G$.  Without loss of
generality we can assume that $D$ is non-separating.  Compress $H$
along $D$ to obtain a new handlebody $H_1$ with boundary $F_1$; let
$G_1$ be the complement of $H_1$.

If $\bdd E$ is inessential in $F_{1}$ then it bounds a disk in $H_{1}
\subset H$ as well, so $F$ is reducible.
       
If $\bdd E$ is essential in $F_{1}$ then $E$ is a compressing disk in
$G_1$ so we can apply the inductive hypothesis to $H_1$.  If 1 or 3
holds then it holds for $H$, and we are done.  Suppose instead $F_1$
is reducible.  Let ${\bf S}$ be a collection of disjoint reducing
spheres for $F_1$ chosen to maximize $\sss$ among all possible such
collections and then, subject to that condition, further choose ${\bf
S}$ to minimize $|E \cap {\bf S}|$.  Clearly $E \cap {\bf S}$ contains
no closed curves, else replacing a subdisk lying in the disk
collection ${\bf S} \cap G_{1}$ with an innermost disk of $E - {\bf S}$
would reduce $|E \cap {\bf S}|$.  Similarly, we have

\bigskip

{\bf Claim 1} Suppose $\eee$ is an arc component of $\bdd E - {\bf S}$
and $F_{0}$ is the component of $F_{1} - {\bf S}$ in which $\eee$
lies.  If $\eee$ separates $F_{0}$ (so the ends of $\eee$ necessarily
lie on the same component of $\bdd F_{0}$) then neither component of
$F_{0} - \eee$ is planar.

{\bf Proof of claim 1:} Let $c_{0}$ be the closed curve component of
$\bdd F_{0} \subset {\bf S} \cap F_{1}$ on which the ends of $\eee$
lie and, of the two arcs into which the ends of $\eee$ divide $c_{0}$,
let $\aaa$ be adjacent to a planar component of $F_{0} - \eee$.  Then
the curve $\eee \cup \aaa$ clearly bounds a disk in both $G_{1}$ and
$H_{1}$ and then so does the curve $c' = \eee \cup (c_{0} - \aaa)$. 
Let $S'$ be a sphere in $N$ intersecting $F_{1}$ in $c'$ and $S_{0}$
be the reducing sphere in ${\bf S}$ containing $c_{0}$.  Replacing
$S_{0}$ with $S'$ (or just deleting $S_{0}$ if $c'$ is inessential in
$F_{1}$) gives a new collection ${\bf S'}$ of disjoint reducing
spheres, intersecting $\bdd E$ in at least two fewer points.  Moreover
$\sss({\bf S'}) = \sss({\bf S})$ since the only change in the
complementary components in $F_{1}$ is to add to one component and
delete from another a planar surface along an arc in the boundary. 
Then the collection ${\bf S'}$ contradicts our initial choice for
${\bf S}$, a contradiction that proves the claim.

\bigskip

Let $H'$ be
the closed complement of ${\bf S}$ in $H_{1}$, so $H'$ is itself a
collection of handlebodies.

\bigskip

{\bf Claim 2} Either $F$ is reducible or $\bdd H'$ is compressible in $N
- H'$.

{\bf Proof of claim 2:} If $\bdd E$ is disjoint from ${\bf S}$ and is
inessential in $\bdd H'$, then $\bdd E$ bounds a disk in $H'$, hence
in $H$, so $F$ is reducible.  If $\bdd E$ is disjoint from ${\bf S}$
and is essential in $\bdd H'$, then $E$ compresses $\bdd H'$ in $N -
H'$, verifying the claim.  Finally, if $E$ intersects ${\bf S}$,
consider an outermost disk $A$ cut off from $E$ by ${\bf S}$. 
According to Claim 1, this disk, together with a subdisk of ${\bf S}$,
constitute a disk $E'$ that compresses $\bdd H'$ in $N - H'$, proving 
the claim.

\bigskip

Following Claim 2, either $F$ is reducible or the inductive
hypothesis applies to a component $H_{0}$ of $H'$.  If 2 holds for
$H_{0}$ then consider a reducing sphere $S$ for $H_{0}$, isotoped so
that the curve $c = S \cap \bdd H_{0}$ is disjoint from the disks
${\bf S} \cap H_{0}$.  The disk $S - H_{0}$ may intersect $H_{1}$; by
general position with respect to the dual 1-handles, each component of
intersection is a disk parallel to a component of ${\bf S} \cap
H_{1}$.  But each such disk can be replaced by the corresponding disk
in ${\bf S} - H_{1}$ so that in the end $c$ also bounds a disk in $N -
H_{1}$.  After this change, $S$ is a reducing sphere for $F_{1}$ in
$N$ and, since $c$ is essential in $H_{0}$, $\sss({\bf S} \cup S) >
\sss({\bf S})$, contradicting our initial choice for ${\bf S}$.  Thus
in fact 1 or 3 holds for $H_{0}$, hence also for $H$.
\end{proof}

In the specific case $N=S^3$, we apply precisely the same argument, 
combined with Waldhausen's theorem \cite{W} on Heegaard splittings of 
$S^3$, to obtain:

\begin{cor}
Let $H$ be a handlebody imbedded in $S^3$, and suppose $G$, the 
complement of $H$, has compressible boundary.  Then either $H$ has an 
unknotted core or the boundary of $H$ is reducible.
\end{cor}

This corollary is similar to (\cite {MT}, Theorem 1.1), but no 
reimbedding of $S^3-H$ is required.

\section{Complements of handlebodies in non-Haken manifolds}

In \cite{F} (see also \cite{MT} for a brief version) Fox showed that
any compact connected $3$-dimensional submanifold $M$ of $S^3$ is
homeomorphic to the complement of a union of handlebodies in $S^3$. 
We generalize this result to non-Haken manifolds, showing that a
submanifold $M$ 
of a non-Haken manifold $N$ has an almost equally simple description,
that is, $M$ is homeomorphic to the complement of handlebodies either
in $S^3$ or in $N$.

\begin{defin}
Let $N$ be a compact irreducible $3$-manifold, and let $M$ be a
compact $3$-submanifold of $N$.  We will say the complement of $M$ in
$N$ is {\em standard} if it is homeomorphic to a collection of
handlebodies or to $N$\#(collection of handlebodies).  (We regard 
$B^{3}$ as a handlebody of genus $0$.)
\end{defin} 

Note that in the latter case $M$ is actually homeomorphic to the
complement of a collection of handlebodies in $S^3$.

\begin{theorem} \label{theorem:complement}
Let $N$ be a closed orientable irreducible non-Haken $3$-manifold, and
let $M$ be a connected compact $3$-submanifold of $N$ with non-empty
boundary.  Then $M$ is homeomorphic to a submanifold of $N$ whose
complement is standard.
\end{theorem}

\begin{proof}
    
The proof will be by induction on $n + g$ where $n$ is the number of
components of $\bdd M$ and $g$ is the genus of $\bdd M$, that is, the
sum of the genera of its components.  If $n + g = 1$ then $\bdd M$ is
a single sphere.  Since $N$ is irreducible, the sphere bounds a
$3$-ball in $N$.  So either $M$ or its complement is a $3$-ball and 
in either case the proof is immediate.

For the inductive step, suppose first that $\bdd M$ has multiple
components $T_{1}, \ldots, T_{n}, n \geq 2$.  Each component $T_{i}$
must bound a distinct component $J_{i}$ of $N - M$ since each must be
separating in the non-Haken manifold $N$.  Let $M' = M \cup J_{n}$; by
inductive assumption $M'$ can be reimbedded so that its complement is
standard.  After the reimbedding, remove $J_{n}$ from $M'$, to recover
a homeomorph of $M$ and adjoin $J_{1}$ (now homeomorphic either to a
handlebody or to $N$\#(handlebody)) instead.  Reimbed the resulting
manifold so that its complement is standard and remove $J_{1}$ to
recover $M$, now with standard complement.
    
Henceforth we can therefore assume that $\bdd M$ is connected and not
a sphere.  Since $N$ is non-Haken there exists a compressing disk $D$
for $\partial M$.

{\bf Case 1.}   $\partial D$ is non-separating on  $\partial M$.

If $D$ lies inside $M$, compress $M$ along $D$ to obtain $M'$  and 
use the induction hypothesis to find an imbedding of $M'$ with 
standard complement.   Reconstruct $M$ by attaching a trivial 
1-handle to $M'$, thus simutaneously attaching a trivial 1-handle to 
the complement.

If $D$ lies outside $M$, attach a 2-handle to $M$ corresponding to $D$
to obtain $M'$, whose connected boundary has lower genus.  Invoking
the inductive hypothesis, imbed $M'$ in $N$ with standard complement. 
Reconstruct $M$ from $M'$ by removing a co-core of the attached
2-handle, thus adding a 1-handle to the complement of $M'$.

{\bf Case 2.}   $\partial D$ is separating on  $\partial M$.

Suppose $D$ lies outside $M$.  Then $D$ also separates $J$ into two
components, $J_1$ and $J_2$, since $H_{2}(N) = 0$.  Denote the
components of $\partial M - \bdd D$ by $\partial_1 \subset J_1$ and
$\partial_2 \subset J_2$, both of positive genus.  Let $M'=M\cup J_2$. 
Reimbed $M'$ so that its complement is standard.  The boundary of $M'$
consists of $\partial_1$ together with a disk.  Since the complement
of $M'$ is standard, there is a non-separating compressing disk $D'$
for $\partial M'$ contained in the complement of $ M'$.  $D'$ is also
a non-separating compressing disk for the reimbedded $\partial M$
(which is contained in $M'$).  Apply case 1 to this new imbedding
of $M$.

We can now suppose that the only compressing disks for $\partial M$
are separating compressing disks lying inside $M$.  Choose a family
$\bf D$ of such $\bdd$-reducing disks for $M$ that is maximal in the
sense that no component of $M' = M - {\bf D}$ is itself
$\bdd$-compressible.  Since each compressing disk is separating,
$genus(\bdd M') = genus(\bdd M) > 0$ so $\bdd M'$ is compressible in
$N$.  Such a compressing disk $E$ can't lie inside $M'$, by
construction, so it lies in the connected manifold $N - M'$; let
$M_{1}$ be the component of $M'$ on whose boundary $\bdd E$ lies. 
Since each disk in $\bf D$ was separating, $M$ has the simple
topological description that it is the boundary-connect sum of the
components of $M'$.  So $M$ can easily be reconstructed from $M'$ in
$N - M'$ by doing boundary connect sum along arcs connecting each
component of $M' - M_{1}$ to $M_{1}$ in $N - (M' \cup E)$.  After this
reimbedding of $M$, $E$ is a compressing disk for $\bdd M$ that lies
outside $M$, so we can conclude the proof via one of the previous
cases.
\end{proof}

\section{Aligned Fox reimbedding}

Now we combine results from the previous two sections and consider
this question: If $M$ is a connected $3$-submanifold of a non-Haken
manifold $N$ and $M$ is $\bdd$-reducible, to what extent can a
reimbedding of $M$, so that its complement is standard, have its
$\bdd$-reducing disks aligned with meridian disks of its complement. 
Obviously non-separating disks in $M$ cannot have boundaries matched
with meridian disks of $N - M$, since $N$ contains no non-separating
surfaces.  But at least in the case when $\bdd M$ is connected, this
is the only restriction.  

\begin{defin} For $M$ a compact irreducible orientable $3$-manifold,
define a disjoint collection of separating $\bdd$-reducing disks ${\bf
D} \subset M$ to be {\em full} if each component of $M - {\bf D}$ is
either a solid torus or is $\bdd$-irreducible. 

For $M$ reducible,
${\bf D} \subset M$ is full if there is a prime decomposition of $M$
so that for each summand $M'$ of $M$ containing boundary, ${\bf D}
\cap M'$ is full in $M'$.

$M \subset N$ a $3$-submanifold is {\em aligned to a standard 
complement} if the complement of $M$ is standard and there is a 
(complete) collection of reducing spheres ${\bf S}$ for $\bdd M$ so 
that ${\bf S} \cap M$ is a full collection of $\bdd$-reducing disks 
for $M$.
\end{defin}

There is a uniqueness theorem, presumably well-known, for 
full collections of disks, which is most easily expressed for 
irreducible manifolds:

\begin{lemma} \label{lemma:natural} Suppose $M$ is an irreducible
orientable $3$-manifold with boundary and $M$ is expressed as a
boundary connect sum in two different ways: $M = M_{1} \natural M_{2}
\natural \ldots \natural M_{n} = M_{1}^{*} \natural M_{2}^{*} \natural
\ldots \natural M_{n^{*}}^{*}$, where each $M_{i}, M_{j}^{*}$ is
either a solid torus or $\bdd$-irreducible.  Then, after
rearrangement, $n^{*} = n$ and $M_{i} \cong M_{i}^{*}$.
\end{lemma}

\begin{proof} One can easily prove the theorem from first principles,
along the lines of e.  g. \cite[Theorem 3.21]{H}, the standard proof
of the corresponding theorem for connected sum.  But a cheap start is
to just double $M$ along its boundary to get a manifold $DM$.  The
decompositions above double to give connected sum decompositions of
$DM$ in which each factor consists of either $S^{1} \times S^{2}$ or
the double of an irreducible, $\bdd$-irreducible manifold which is
then necessarily irreducible.  Then \cite[Theorem 3.21]{H} implies that
$n = n^{*}$ and that the two original decompositions of $M$ also each
contain the same number of solid tori.  After removing these, we are
reduced to the case in which the only $\bdd$-reducing disks in $M$ are
separating and $n^{*} = n$.

Following the outline suggested by the proof of \cite[Theorem
3.21]{H}, choose a disk $D$ that separates $M$ into the component
$M_{n}$ and the component $M_{1} \natural M_{2} \natural \ldots \natural
M_{n-1}$.  Choose disks $E_{1}, \ldots, E_{n-1}$ that separate $M$
into the components $M_{1}^{*} , M_{2}^{*} , \ldots M_{n^{*}}^{*}$. 
Choose the disks to minimize the number of intersection components in
$D \cap (\cup \{ E_{i} \})$.  Since each manifold is irreducible and
$\bdd$-irreducible, a standard innermost disk, outermost arc argument
(in $D$) shows that in fact $D$ is then disjoint from $\{ E_{i} \}$,
so $D \subset M_{n}^{*}$ (say).  Since $M_{n}^{*}$ is
$\bdd$-irreducible, $D$ is $\bdd$-parallel in $M_{n}^{*}$ so in fact
(with no loss of generality) $M_{n} \cong M_{n}^{*}$ and $M_{1}
\natural M_{2} \natural \ldots \natural M_{n-1} \cong M_{1}^{*} \natural
M_{2}^{*} \natural \ldots \natural M_{n-1}^{*}$.  The result follows by
induction.  \end{proof}

\begin{theorem} \label{theorem:align}
 Let $N$ be a closed orientable irreducible non-Haken $3$-manifold,
 and $M$ be a connected compact $3$-submanifold of $N$ with connected
boundary.  Then $M$ can be reimbedded in $N$ with standard
 complement so that $M$ is aligned to the standard complement.
\end{theorem}

\begin{proof} The proof is by induction on the genus of $\bdd M$. 
Unless $M$ has a separating $\bdd$-reducing disk, there is nothing
beyond the result of Theorem \ref{theorem:complement} to prove. So we
assume that $M$ does have a separating $\bdd$-reducing disk; in
particular the genus of $\bdd M$ is $g \geq 2$.  We inductively assume
that the theorem has been proven whenever the genus of $\bdd M$ is
less than $g$.

The first observation is that it suffices to find an embedding of $M$
in $N$ so that there is some reducing sphere $S$ for $\bdd M$ in $N$. 
For such a reducing sphere divides $J = N - M$ into two components
$J_{1}$ and $J_{2}$.  Apply the inductive hypothesis to $M \cup J_{1}$
to reimbed it with aligned complement $J_{2}'$.  Notice that by a
standard innermost disk argument, the reducing spheres can be taken to
be disjoint from $S$.  After this reimbedding, apply the inductive
hypothesis to $M \cup J_{2}'$ to reimbed it so that its complement
$J_{1}'$ is aligned.  After this reimbedding, $M$ has aligned
complement $J_{1}' \cup_{S - M} J_{2}'$.

Our goal then is to find a reimbedding of $M$ so that afterwards $\bdd
M$ has a reducing sphere.  First use Theorem \ref {theorem:complement}
to reimbed $M$ in $N$ so that its complement $J$ is standard, i.  e.
either a handlebody or $N$\# (handlebody).  Since $M$ is
$\bdd$-reducible, Lemma \ref{lemma:pairodisks} applies: either $M$ is
itself a handlebody (in which case the required reimbedding of $M$ is
easy) or there are disjoint compressing disks $D$ in $J$ and $E$ in
$M$.  Since $J$ is standard, $D$ can be chosen to be non-separating in
$J$.  Then $\bdd E$ is not homologous to $\bdd D$ in $\bdd M$ so $\bdd
E$ is either separating in $\bdd M$ or non-separating in $\bdd M -
\bdd D$.  In the latter case, two copies of $E$ can be banded together
along an arc in $\bdd M - \bdd D$ to create a separating essential
disk in $M$ that is disjoint from $D$.  The upshot is that we may as
well assume that $D \subset J$ is non-separaring and $E \subset M$ is
separating.

Add a $2$-handle to $M$ along $D$ to get $M'$, still with standard
complement $J'$.  Dually, $M$ can be viewed as the complement of the
neighborhood of an arc $\aaa \subset M'$.  If $\bdd E$ is inessential
in $\bdd M'$, it bounds a disk $D'$ in $J' \subset J$.  Then the
sphere $D' \cup E$ is a reducing sphere for $M$ as required.  So we
may as well assume that $\bdd E$ is essential in $\bdd M'$ and of
course still separates $M'$.  By inductive assumption $M'$ can be
embedded in $N$ so that its complement is aligned, but note that this
does not immediately mean that $\bdd E$ itself bounds a disk in $N -
M'$.  Let ${\bf S}$ be a complete collection of reducing spheres for
$\bdd M'$ intersecting $M'$ in a full collection of disks.

$E$ divides $M'$ into two components, $U$ and $V$ with, say, $\aaa
\subset U$.  If $M'$ is reducible (i.e. contains a punctured copy of
$N$) an innermost (in $E$) disk argument ensures that the reducing
sphere is disjoint from $E$.  By possibly tubing $E$ to that reducing
sphere, we can ensure that the $N$-summand, if it lies in $M'$, lies
in $U \subset M'$.  That is, we can arrange that $V$ is irreducible. 
$E$ extends to a full collection of disks in $M'$, with the new disks
dividing $U$ and $V$ into $\bdd$-connected sums: $U = U_{1} \natural
\ldots \natural U_{m}, V = V_{1} \natural \ldots \natural V_{n}, m, n
\geq 1$, with each $U_{i}, V_{j}$ either $\bdd$-irreducible or a solid
torus (with one of the $U_{i}$ possibly containing $N$ as a connect
summand).  By Lemma \ref{lemma:natural}, some component $V'$ of $M' -
{\bf S}$ is homeomorphic to $V_{n}$.  Tube together all components of
${\bf S}$ incident to $V'$ along arcs in $\bdd V'$ to get a reducing
sphere $S'$ dividing $M'$ into two components, one homeomorphic to
$V_{n}$ and the other homeomorphic to $U \natural V_{1} \natural V_{2}
\natural \ldots \natural V_{n-1}$.  The latter homeomorphism carries
$\aaa \subset U$ to an arc $\aaa'$ that is disjoint from the reducing
sphere $S'$.  Then $M' - \eta(\aaa')$ is homeomorphic to $M$ and
admits the reducing sphere $S'$.  In other words, the reimbedding of
$M$ that replaces $M' - \eta(\aaa)$ with $M' - \eta(\aaa')$ makes $\bdd
M$ reducible in $N$, completing the argument.  \end{proof}

\begin{cor} Given $M \subset N$ as in Theorem \ref{theorem:align}, 
suppose $\bf D$ is a full set of disks in $M$.  Then, with at most one 
exception, each component of $M - \bf D$ embeds in $S^{3}$.
\end{cor}

\begin{proof} Following Theorem \ref{theorem:align} reimbed $M$ in $N$
with standard complement so that $M$ is aligned to the standard
complement.  Then there is a collection $\bf S$ of disjoint spheres in
$M$ so that, via Lemma \ref{lemma:natural}, $M - {\bf S}$ and $M - 
{\bf D}$ are homeomorphic.  Since $N$ is irreducible, each 
component but at most one of $N - {\bf S}$ is a punctured $3$-ball.  
Finally, each component of $N - {\bf S}$ contains at most one 
component of $M - {\bf D}$ since each component of $\bf S$ is 
separating. \end{proof}

\end{document}